\documentclass[leqno,12pt]{article}

\usepackage{graphicx}\usepackage{epsfig,graphics}
\usepackage[scale=0.8]{geometry}
\usepackage[T1]{fontenc}
\usepackage{svg}
\usepackage{comment}
\usepackage{hyperref}
\usepackage{fancyhdr}
\usepackage{lastpage}
\usepackage{algorithm,algpseudocode}
\usepackage{theorem}
\usepackage{amsmath}
\usepackage{multirow}

\def\epsilon{\varepsilon} \def\phi{\varphi}       \def\reals{{I\kern-.35em R}}  
\def\argmax{\mathop{\rm argmax}}  
\def\lset{\big\{\,}    \def\mset{\,\big|\,}   \def\rset{\,\big\}}
\def\Lset{\Big\{\,}    \def\Mset{\,\Big|\,}   \def\Rset{\,\Big\}}
\outer\def\proclaim #1. #2\par{\medbreak \noindent{\bf#1.\enspace}{\sl#2}\par 
   \ifdim\lastskip<\medskipamount \removelastskip\penalty55\medskip\fi} 
\def\downto{{\raise 1pt \hbox{$\scriptstyle \,\searrow\,$}}} 
\def\upto{{\raise 1pt \hbox{$\scriptstyle \,\nearrow\,$}}} 
\def\text#1{\;\,\hbox{#1}\;\,}  \def\txt#1{\,\hbox{#1}\,}  
\def\state #1. { \noindent{\bf#1.\enspace}} 
 
\def\eqalign#1{\begin{array}{lcr} #1 \end{array}} 
\def\low#1{{\lower1pt \hbox{$\scriptstyle #1$}}}   
\def\comp{\raise 1pt \hbox{$\scriptstyle\circ$}} 
\def\half{{{}\raise 1pt\hbox{${{\scriptstyle1}\over{\scriptstyle 2}}$}}} 
     \def\sumn{\sum\nolimits}
\def\implies{\quad\hbox{$\Longrightarrow$}\quad} 
\def\tto{\;{\lower 1pt \hbox{$\rightarrow$}}\kern -12pt
           \hbox{\raise 2.8pt \hbox{$\rightarrow$}}\;}
\def\substack#1#2{{\scriptstyle{#1}\atop\scriptstyle{#2}}} 
  
\def\mdot{{\kern-.02em\cdot\kern-.04em}} 
\def\newline{\vskip0cm\noindent}

\def\pls{{\scriptscriptstyle +}} \def\mns{{\scriptscriptstyle -}}
\def\eop{\hfill{$\vcenter{\hrule height1pt \hbox{\vrule width1pt height5pt  
   \kern5pt \vrule width1pt} \hrule height1pt}$} \medskip} 
 
\def\cC{{\cal C}}

\begin{document}
\centerline{\large \bf REACHING AN EQUILIBRIUM OF PRICES AND HOLDINGS}
\smallskip
\centerline{\large \bf OF GOODS THROUGH DIRECT BUYING AND SELLING } 
\bigskip
\bigskip
\centerline{\bf J.\ Deride }
\centerline{Universidad T\'ecnica Federico Santa Maria}
\centerline{Santiago, Chile   
(julio.deride@usm.cl)}
     \medskip
\centerline{\bf A.\ Jofr\'e }
\centerline{Center for Mathematical Modeling, Univ.\ of Chile}
\centerline{Casilla 170/3, Correo 3, Santiago, Chile  
(ajofre@dim.uchile.cl)}
     \medskip
\centerline{\bf R.\ T.\ Rockafellar }
\centerline{Dept. of Mathematics, Univ.\ of Washington, Seattle, WA 98195-4350
(rtr@uw.edu)}

\bigskip
\bigskip
\bigskip
\state Abstract.
The Walras approach to equilibrium focuses on the existence of market 
prices at which the total demands for goods are matched by the total 
supplies.  Trading activities that might identify such prices by bringing 
agents together as potential buyers and sellers of a good are 
characteristically absent, however.  Anyway, there is no money to pass from 
one to the other as ordinarily envisioned in buying and selling.  Here a 
different approach to equilibrium --- what it should mean and how it may be 
achieved --- is offered as a constructive alternative.

Agents operate in an economic environment where adjustments to holdings
have been needed in the past, will be needed again in a changed future,
and money is familiar for its role in facilitating that.  Marginal utility 
provides relative values of goods for guidance in making incremental 
adjustments, and with money incorporated into utility and taken as 
num\'eraire, those values give money price thresholds at which an agent 
will be willing to buy or sell.  Agents in pairs can then look at such 
individualized thresholds to see whether a trade of some amount of a good 
for some amount of money may be mutually advantageous in leading to higher 
levels of utility.  Iterative bilateral trades in this most basic sense, 
if they keep bringing all goods and agents into play, are guaranteed in 
the limit to reach an equilibrium state in which the agents all agree on 
prices and, under those prices, have no interest in further adjusting their 
holdings.  The results of computer simulations are provided to illustrate 
how this works.


\bigskip
\bigskip
\centerline{Version of 23 January 2023}

\newpage

\bigskip

\section{
               Introduction}

For agents with holdings in various goods and preferences given by utility
functions, exchanges in quantities of the goods may be undertaken to mutual 
advantage. Those exchanges might be facilitated by a ``market,'' but how 
would that operate and what would it achieve?  

In the usual view of a market, there are prices at which an agent can
buy or sell goods subject to the budget constraint that the total value of
the purchases is covered by the total value of the sales.  An {\it 
equilibrium of prices and holdings\/} (in the terminology we use here) is 
at hand when no agent has a prospect of buying and selling that could lead 
to holdings with a higher level of utility.  What might bring this about?  
Could agents proceed directly with each other, making advantageous exchanges 
and, in so doing, gradually reach a consensus on prices and thereby an 
equilibrium?

Answering this question seems essential for understanding the very 
fundamental concept of a market, but past research on equilibrium has 
mostly taken a different path.  The standard line of 
inquiry, opened by Walras in the 1870s and placed on firm mathematical
footing by Debreu \cite{Debreu-Value} in 1959, is whether, for the initial 
holdings of the agents, there exist prices such that, when each agent 
individually determines an exchange that would maximize utility under the 
budget constraint, the resulting aggregate demands for goods would match 
the overall supplies.  Then indeed, after those desired adjustments, the 
prices and holdings will be in equilibrium, with no agent interested in 
further exchange.  A balance in competition is depicted this way as arising 
from individual agents acting solely out of self-interest.  Prices  
coordinate them and appear to serve as a means of decentralization.  But 
there are flaws in this portrayal.

Although the existence of Walrasian prices has been established under 
various mathematical assumptions, explaining how they might actually come 
into being has been harder.  The main rationale has been some version of 
t\^atonnement, in which prices are adjusted up or down by an abstract entity 
who repeatedly seeks feedback from the agents.  This process is usually
modeled by a differential equation;  see 
\cite{Uzawa, Balasko0, Saari, Walker, Keisler, Kitti}, and more recently
\cite{Stability}.   In certain cases t\^atonnement does identify prices 
fitting the Walras prescription,\footnote{ 
    Contrary to wide-spread pessimism about the ability of the process to
    succeed, it can generally be counted on to converge to Walrasian prices 
    when initiated from holdings not too far from equilibrium holdings,
    even in a setting where agents are allowed to hold zero quantities of
    some of the goods; see \cite{Stability}. }
yet it suffers from requiring a central coordinator instead of relying on 
agent interactions directly.  And it lacks justification as an information 
process backed by economic observation.  

Another serious shortcoming of the Walrasian market is that the picture of 
buying and selling is left blank.  When people think of a real market, they 
contemplate circumstances in which someone can buy or sell a quantity of 
some good {\it for a quantity of money}, but in the Walrasian tradition no 
money changes hands.  All talk of money is in fact set aside.  Prices are 
only {\it relative}, with their ratios indicating the values of goods 
relative to each other.   The agents, having determined what they want to 
acquire or give up, while respecting the budget constraint (which is 
insensitive to the scaling behind the relative prices), present their 
decisions to a ``clearing house entity.''  That entity then somehow executes 
the desired transfers of goods as a simultaneous {\it grand exchange}.  This 
is not at all a matter of direct buying or selling by agents, and it can 
hardly be viewed as signaling economic decentralization.   

Might support for the Walras approach be found by demonstrating that the 
grand exchange could be carried out by the agents themselves, trading 
together directly in some order?  An investigation by one of this paper's 
coauthors (Rockafellar) in \cite{WalrasExchanges} showed that this might
be possible up to a point, with two agents at a time bartering one good for
another at the given prices.  But there was a hitch:  no guarantee that
each trade would result in higher utility for both.  With agents permitted
to barter {\it bundles of several\/} goods in each transaction, a 
utility-raising trade was determined in \cite{WalrasExchanges} always to 
exist, but whether a sequence of such trades could bring about the 
specified grand exchange was elusive and left unanswered.  In fact, this 
was only the latest inquiry into a genuine market justification for 
Walrasian prices.  Starr \cite{Starr} (1972) and Ostroy \cite{Ostroy} (1973) 
already derived conditions under which the grand exchange might be achieved 
through budget-balanced multigood trades among agents with money in a 
helping role.  Their trades, though, don't meet the criterion underscored 
in \cite{WalrasExchanges} that utility should always go up for both parties.

Doubts about the credibility of the Walras scheme led others to abandon it 
in favor of alternative models of exchange and price formation.  Shapley 
and Shubik \cite{ShapleyShubik} in 1977 had agents make money-bids for 
quantities of goods which trading posts for those goods would settle at 
auction-fixed prices.  This replaced the centralized clearing house by a 
noncooperative game with utility pay-offs.   It marked the beginning an 
elaborate theory of {\it strategic market games\/}, envisioned as perhaps 
building a bridge from microeconomics to macroeconomics.   For a bigger 
view of this over the years, see Dubey and Shubik \cite{DubeyShubik} (1980), 
Sorin \cite{Sorin} (1996), and Levando \cite{Levando} (2012).  A strategic 
market game achieves a different equilibrium than the one in the Walras 
model, and that equilibrium might not even provide an allocation of goods 
that is ``efficient'' in the sense of Pareto optimality with respect to 
agent utility.  Getting to equilibrium would seem moreover to require a 
degree of organization and rule enforcement that raises many 
questions.\footnote{
   For instance, agents in \cite{ShapleyShubik} must back up their bids
   with money commitments before knowing what the equilibrium prices might 
   turn out to be, and this can present issues of insolvency that have to 
   be handled.  }
It's hard not to see a plausiblity gap with respect to what really goes on 
in economic behavior.  

Still other researchers have departed from the Walras model and the prices 
it offers without passing to game formulations.  How might agents, through 
repeated bilateral exchanges to the mutual benefit of both, eventually 
achieve Pareto optimality and with it some particular equilibrium of prices 
and holdings?  This is the direction that especially interests us here.

Feldman \cite{Feldman} in 1973 addressed this with each pair of agents 
persisting in multigood trades until all possibilities of improvement have 
been exhausted.  Only then could another pair take the stage.  He showed 
that cluster points of the sequence of goods allocations would be Pareto 
optimal as long as one of the goods was ``money-like'' in a sense.  
Eckalbar \cite{Eckalbar} in 1986 improved the picture greatly by focusing 
on pairs of agents who {\it exchange a single good at a time for money} --- 
with one being the buyer and the other the seller.   The marginal utility 
of money relative to the marginal utility of the good, for each of the 
agents, guides whether the deal is worthwhile.  This is the pattern we will 
adopt as well, in our developments below, but our mechanism for reaching an 
equilibrium will be very different from Eckalbar's, which was a process in 
continuous time rather than one of iterative discrete exchange.  He poses 
a differential equation in which the right side is supposed to have certain 
behaviors (not fully pinned down in terms of the properties of the function 
itself), and the claim is made that trajectories generated from this will 
converge to Pareto optimality.  However, his proof of this claim falls 
short.\footnote{
     It takes for granted that the trajectory, being bounded, must converge
     to a particular point, whereas it conceivably might just
     approach a ``limit cycle.''}
Another approach to reaching an equilibrium through continuous trading 
governed by a differential equation was proposed by Bottazzi \cite{Bottazzi}
in 1994.  There, trading is bilateral in a generalized sense, but it 
doesn't come down to simple buying/selling of individual goods and indeed, 
there is no money.  A drawback to both of these differential equation 
schemes, from the perspective of economics, is that they tacitly require 
some societal entity to select the particular equation and compel the 
agents to submit to it.  Bottazzi also expects agents to wish for 
``steepest ascent'' with respect to utility, but that is highly artificial,
because the direction of steepest ascent varies with the choice of units in 
which the goods are measured.  Why should the choice of units affect an
agent's behavior?

Bilateral trading in discrete time has been taken up more recently by 
Fl{\aa}m in \cite{Flaam1}, \cite{Flaam2}, \cite{Flaam3}.\footnote{
    Coauthor Rockafellar of the current paper is grateful to Fl{\aa}m
    for introducing him to this subject during a two-week visit to the 
    University of Bergen in June of 2015.} 
In \cite{Flaam1} utility is somehow measured directly in money units 
and trading is accompanied by monetary side payments.  Like Bottazzi, 
Fl{\aa}m insists on directions of ``steepest ascent'' for trades between 
two agents.\footnote{
    He relies on the Euclidean norm for this.  But that amounts to 
    identifying the space of goods vectors with the dual space of price 
    vectors.  This has no economic basis and thus undermines the economic 
    rationale for his trading prescriptions.}
Besides the artificiality of those directions, they generally force agents 
to trade ``full bundles'' of goods, instead of a single good at a time.   
Bilateral trading of a different sort, not involving money, was taken up 
from another angle by Fl{\aa}m and Gramstad in \cite{Flaam2}.  
A subsequent contribution of Fl{\aa}m in \cite{Flaam3} brings in money as 
num\'eraire and proposes a scheme in which pairs of agents trade bundles of 
goods at money prices, but without being obliged to exchange money itself.   
Negative prices and negative holdings of money are permitted, and besides 
having utility functions, pairs of agents are viewed as possessing various 
{\it joint\/} properties to guide their trading with each other that aren't 
normally a part of equilibrium models.  Despite those innovations, the 
assumptions and arguments in \cite{Flaam3} don't rise to the level of 
proving that the iterative adjustments will converge to an equilibrium 
configuration, or for that matter even produce an equilibrium as a cluster 
point.  Another thing to note is that a transaction between two agents in 
\cite{Flaam3} isn't required to be advantageous to both parties.   One of 
the agents could actually see a decrease in utility, and the question then 
comes up as to why the transaction would be carried out.  It seems that an 
enforcer would be needed in the background.  

Most recently in \cite{Flaam4, Flaam5}, Fl{\aa}m has explored different 
money-price schemes in a game framework based once more on a Walras-like 
clearing house operation.  Double auctions are employed in a process of 
identifying prices that might ultimately signal Pareto optimality, but 
bilateral trades are no longer seen.  This draws in part on ideas of Gintis 
and Mandel \cite{Gintis, Mandel}, who have agents participating in a game 
where they choose price vectors as strategies and get utility pay-offs 
from the resulting redistributions of goods coming from processes in a 
centralized clearing house system.  More about double auctions in 
determining a price for a good, and the stability of that process or lack 
of it, can be learned from current work of Rasooly \cite{Rasooly}.

In contrast, our approach in this paper is far simpler, yet demonstrates 
how an equilibrium of prices and holdings is sure to emerge in time from 
elementary bilateral transactions which agents can undertake solely for 
reasons of immediate self-interest.  Another major difference is that, 
instead of multigood trades, we emphasize like Eckalbar the {\it buying and 
selling of a single good at a time} and insist on {\it a payment of money 
always being made}.  The amount of payment emerges from utility of the two 
agents involved.  We are able to show that interactions on this utterly 
fundamental level lead inevitably to prices accepted by all the agents,
which identify a particular equilibrium of prices and holdings --- and
{\it as a true limit}, not merely some cluster point.  Pareto optimality 
of the agents' ultimate holdings is revealed in this way to be a natural 
consequence of direct deals that have no need for a clearing house 
entity.\footnote{
    Because an equilibrium of prices and holdings can be viewed as the 
    case of a Walras equilibrium in which the initial holdings coincide 
    with the final holdings, its Pareto optimality is assured by the First 
    Welfare Theorem.}
Nothing like that has been established before now.\footnote{    
    These results, written up in a working paper of April 2018, were 
    presented by co-author Jofr\'e in June 2018 at a meeting of economists 
    in Paris.}

Another distinctive feature of our contribution is that we accompany the 
money-based trading between agents by computational experiments which 
illustrate how it works out.   In a numerical model with Cobb-Douglas-type 
utility functions, we show that an exchange equilibrium can readily be 
attained through the decentralized process we propose, even with many goods 
and many agents.  Different trading equilibria are seen to be reached by 
the agents even in starting from the same holdings, because of the random 
patterns in which the agents can come together in pairs and select a good 
to buy/sell.  This is completely ``economic'' and underscores that our 
equilibrium prices and holdings come from decentralized market dynamics 
without a higher entity dictating a plan rigid about who-does-what-when.  
The Walras approach, in contrast, lacks representation of mutually 
beneficial exchange activity between agents in pairs.  Very differently, 
the prices it assigns to a configuration of initial holdings rely on having 
a centralized mechanism for an all-at-once redistribution of goods, as 
explained earlier.  That makes Walrasian prices economically unconvincing,
and not just to us, as the literature reveals.  The realization that simple 
buying and selling, as we describe it, doesn't lead to Walrasian prices, 
unless by accident, highlights the artificiality of those prices all the more.
\medskip

In the plan for the rest of this paper, we start in Section 2 by explaining 
the mathematical details of our model of utility and what extra help it 
offers to the agents.  We continue in Section 3 by showing how agents are 
able to identify mutually beneficial opportunities for acting as buyer or 
seller of a good, and how the absence of such opportunities characterizes 
a state of equilibrium.  Section 4 develops how iterative buying and 
selling produces convergence to such equilibrium.  This is then brought to
life by computer simulation in Section 5.

\section{
     Marginal utility and price thresholds}

The goods space in our setting is the nonnegative orthant 
$\reals_\pls^{n+1}$.  The goods are indexed by $j=0,1,\ldots,n$, with 
good 0 standing for money.  The agents, indexed by $i=1,\ldots,m$, deal 
with holdings of goods represented by vectors 
$x_i=(x_{i0},x_{i1},\ldots,x_{in})$, and they have utility functions 
$u_i$ for comparing such holdings.  Money will serve as num\'eraire, but 
it will have a role for us beyond that of a customary num\'eraire good
like gold, because the buyer of any other good will have to transfer to 
the seller a quantity of money in payment.

Although it may be deemed controversial in some quarters to have utility 
apply to money, this notion goes back to the early days of the topic we are 
addressing.  We have argued moreover, in previous papers devoted to economic 
equilibrium \cite{GEI,Stability}, that there are solid reasons for allowing 
it, citing strong support from Keynes.  Much of the resistance to money is 
likely due anyway to the Walrasian model of equilibrium as a circumstance 
with no past and no future.  In thinking of it instead as a transient 
phenomenon in which agents agree on prices, but no one wishes to buy or sell 
anything, the picture is different.  Equilibrium holdings can be affected 
by outside actions like consumption, production, taxation or subsidy.  The 
resulting disequilibrium that can then be countered by renewed buying 
and selling.  From that angle, agents can well have an understanding of 
money and its ongoing importance in transactions, and may therefore wish 
to hold a quantity of it.   The 2018 paper of Dubey, Sahi and Shubik
\cite{DubeySahiShubik} provides additional insights into money along with
other helpful references.

Utility functions will do more for us than just capturing the preferences 
that an agent might have for one vector of holdings over another.   They 
will be important also in assessing quantitatively the effects of 
{\it gradual\/} shifts in holdings.  That will give us a handle on
how trading might be executed when two agents get together.  
Specifically, suppose agent $i$ has current holdings given by $x_i$ 
and is contemplating a shift from that to $x_i+\tau\Delta x_i$ for some 
``step size'' $\tau\geq 0$ with respect to a vector $\Delta x_i$ having 
components that may be positive, negative, or zero.   The choice of $\tau$ 
must take into account the utility value that corresponds to it, namely
$$
      \theta(\tau) = u_i(x_i+\tau\Delta x_i).
\eqno(2.1)$$
For our purposes, we want not only to know whether a choice $\tau_1$ is
better than a choice $\tau_0$, in the sense that 
$\theta(\tau_1)>\theta(\tau_0)$, but also to be able to work with derivatives 
$$
     \theta'(\tau) = \nabla u_i(x_i+\tau\Delta x_i)\mdot \Delta x_i
\eqno(2.2)$$ 
as capturing {\it marginal utility\/} with respect to an increase or decrease 
in $\tau$.  

How realistic is that, though, from the perspective of the basic theory of 
preference relations and their utility representations?  Standard axioms only 
produce a utility function that is continuous, quasi-concave and perhaps 
quasi-smooth in the sense that its convex upper level sets have kink-free 
boundaries.  Differentiability is unaddressed, and derivatives might not be 
quantitatively meaningful anyway because of the potential arbitrariness in 
nonlinear rescaling of utility values.  However, recent advances in 
\cite{Preferences} are illuminating in this respect.  Natural properties 
{\it of the preference relation itself\/} are identified there which 
characterize whether it can be represented by a utility function that is 
$\cC^1$ or $\cC^2$.  In the $\cC^2$ case, if the preference sets are 
{\it strongly\/} convex (strictly convex with reliable curvature), there is 
sure to exist, relative any compact convex subset $C$ of the positive 
orthant, a $\cC^2$ utility representation that is concave on $C$, moreover 
in a minimal way that makes it {\it unique up to affine rescaling}, i.e., 
up to the choice of units in which utility is to be measured.  We align 
ourselves here with those preference-relation-grounded characteristics.

In preparation for stating the assumptions, we assign to each agent $i$ a 
set $X_i$ of {\it admissible\/} holdings $x_i$, taking it to have the form
$$
  X_i= \lset x_i\in\reals^{n+1}_\pls \mset x_{ij}>0 \,\txt{for} j\in J_i\rset 
   \text{for a collection $J_i$ of goods $j$.}   
\eqno(2.3)$$
The goods in $J_i$ will be called the {\it essential\/} goods for agent $i$, 
with other goods being {\it inessential}.  Money is assumed to always be 
essential:
$$ 
      0\in J_i \text{for all agents $i$.}
\eqno(2.4)$$

\proclaim Utility assumptions.
Agent $i$ has on $X_i$ a concave utility function $u_i$ that is twice 
continuously differentiable.\footnote{
    The interpretation of this for points of $X_i$ on its boundary is that
    the first and second derivatives on the interior of $X_i$ extend to them
    as unique limits.}
The gradient vectors $\nabla u_i(x_i)$ have positive components and the 
hessian matrix $\nabla^2 u_i(x_i)$ is negative-definite on the subspace 
orthogonal to $\nabla u_i(x_i)$.  Furthermore, 
$$
   x_i\in X_i \implies \lset x'_i \in X_i \mset u(x'_i)\geq u(x_i)\rset
               \text{is a closed subset of $\reals^{n+1}$.}
\eqno(2.5)$$

The condition assumed on gradients guarantees that $u_i$ increases on 
$X_i$ with respect to any increase in  good:
$$
   u_i(x_i') > u_i(x_i)  \text{when} x_i' \geq x_i,\; x_i' \neq x_i.
\eqno(2.6)$$
The concavity and differentiability assumption makes the expressions 
$\theta(\tau)=u_i(x_i+\tau\Delta x_i)$ in (2.1) be concave functions of
$\tau$ with continuous second derivatives, 
$$
   \theta''(\tau)= 
   \Delta x_i\mdot\nabla^2 u_i(x_i+\tau\Delta x_i)\Delta x_i.
$$
The partial negative-definiteness condition on hessians\footnote{ 
     This corresponds to locally strong convexity of upper level sets 
     of $u_i$, instead of the more commonly assumed strict convexity;
     see \cite[Section 4]{Preferences}.}   
assures that $\theta''(\tau)<0$ unless $\Delta x_i$ is a multiple of 
$\nabla u_i(x_i+\tau\Delta x_i)$, 
which will not come up because we will
only be looking at adjustment vectors $\Delta x_i$ that have both a positive
component and a negative component.  Thus, the functions $\theta(\tau)$
we work with will always be {\it strongly\/} concave.  The closedness in 
(2.5) ensures further that, for $x_i\in X_i$ and mixed-sign vectors
$\Delta x_i$, the set
$$
   \lset \tau\geq 0 \mset x_i +\tau\Delta x_i \in X_i,\;
    u_i(x_i +\tau\Delta x_i)\geq u_i(x_i)\rset \text{is compact.}
\eqno(2.7)$$

A key role will be played by the partial derivative 
$\frac{\partial u_i}{\partial x_{i0}}(x_i)$, which gives the 
{\it marginal utility of money\/} associated with the holdings vector 
$x_i$.  It is the derivative $\theta'(0)$ in (2.1) in the case of $\Delta x_i$
having $1$ as its initial component, but $0$ for all the others.  More 
generally, $\frac{\partial u_i}{\partial x_{ij}}(x_i)$ gives the 
{\it marginal utility of good $j$\/} at $x_i$.  These marginal utilities 
are always positive under our assumptions.

\proclaim Price thresholds for goods.  For an agent $i$ and a good $j\neq 0$, 
the ratio
$$
         p_{ij}(x_i) = \frac{\partial u_i}{\partial x_{ij}}(x_i)\Bigg/\
              \frac{\partial u_i}{\partial x_{i0}}(x_i) >0
\eqno(2.8)$$
will be called the price threshold of agent $i$ for selling or buying 
good $j$ when the current holdings are given by the goods vector $x_i$.

This terminology will be justified in the analysis of selling and buying
that comes next.  Observe that when the the marginal utility 
$\frac{\partial u_i}{\partial x_{ij}}(x_i)$ is interpreted as
being measured in units of utility of agent $i$ per unit of good $j$, 
and $\frac{\partial u_i}{\partial x_{ij}}(x_i)$ is interpreted as being 
measured in units of  utility of agent $i$ per unit of money, the prices 
$p_{ij}(x_i)$ correctly come out in units of money per unit of 
good $j$.\footnote{
   Note that the monotonicity of the partial derivatives of $u_i$, following
   from our assumption that $u_i$ increases with respect to upward shifts of 
   components of $x_i$, doesn't necessarily carry over to monotonicity of
   the prices $p_{ij}(x_i)$, due to the division in (2.8). }
It's worth noting that {\it the threshold values $p_{ij}(x)$ are inherent
in the underlying preference relation and independent of the particular
scaling of the function $u_i$};  rescaling can change the length of gradient
vectors but has no effect on the ratios between the components of those
vectors.

We postulate, as is usual in models of exchange aimed at equilibrium analysis,
that agents have in $X_i$ initial holdings $x_i^0$ which they may wish to 
trade for other holdings they would like better.   All goods are assumed to 
be in positive {\it total\/} supply, and those amounts won't be altered by 
the trading.  The trades will be executed in such a way that the utility of 
an agent's holdings never decreases.  Because of (2.7), the holdings of
agent $i$ will always be in $X_i$.  Then in particular, through (2.4) and
(2.5), agents will always have some money at their disposal.  

This feature leads to a major difference between our model and previous
models.  We not only approach trading iteratively through pairs of agents, 
but insist on {\it bringing money into every transaction\/}.  The 
considerations behind that will be laid out now in terms of agents selling
or buying an amount of some good $j\neq 0$ for some amount of the money
good $j=0$.
\medskip

{\bf Agents as sellers.}  For an agent $i$ and a good $j\neq 0$, what
needs to be considered in contemplating the sale of a quantity $\xi_j>0$ 
of that good at a price $\pi_j>0$?  This consists of subtracting $\xi_j$
from the amount of good $j$ currently held, and therefore must be subject
to the constraint $\xi_j \leq x_{ij}$.  On the other hand, it involves
adding the money amount $\pi_j\xi_j$ to the current holding of money,
namely $x_{i0}$.  Thus, it shifts $x_i$ to $x_i +\xi_j[\pi_j,-1]$, where
the notation is that
$$
 [\pi_j,-1] \;=\! \text{the vector $\Delta x_i$ having components} 
   \left\{\eqalign{ \Delta x_{i0} = \pi_j, &\cr \Delta x_{ij} = -1, &\cr
       \Delta x_{ik} = 0 \text{for other goods} k. &} \right.  
\eqno(2.9)$$
In evaluating the utility of such a sale, we are in the framework of (2.1)
for this special $\Delta x_i$ with $\xi_j$ in place of $\tau$ and are
looking at
$$
     \theta_\pls(\xi_j)= u_i(x_i +\xi_j[\pi_j,-1]\,).
\eqno(2.10)$$

For selling to be attractive, starting from an arbitrarily small amount,
the marginal utility $\theta'_\pls(0)$ should be $>0$.  But
$$
  \theta'_\pls(\xi_j) = 
        \frac{\partial u_i}{\partial x_{i0}}(x_i +\xi_j[\pi_j,-1]\,)\pi_j
            - \frac{\partial u_i}{\partial x_{ij}}(x_i +\xi_j[\pi_j,-1]\,),
\eqno(2.11)$$
so that
$$
   \theta'\pls(0) = 
     \frac{\partial u_i}{\partial x_{i0}}(x_i)\Big[\pi_j -p_{ij}(x_i)\Big],
\eqno(2.12)$$ 
where the initial factor, the marginal value of money, is positive.  Therefore, 
$$\eqalign{
\text{selling good $j$ at price $\pi_j$ will be attractive to agent}&\cr 
\text{$i$ with current holdings $x_i$ if and only if 
     $\,\pi_j > p_{ij}(x_i)$.} &}
\eqno(2.13)$$

When selling is attractive in this way, how much should agent $i$ consider
selling?   This is where our assumptions on second derivatives come in.
Under those assumptions the function $\theta(\xi_j)$ in (2.6) is concave
with $\theta'_\pls(\xi_j)$ strictly decreasing from its initial value
$\theta'_\pls(0)>0$ as $\xi_j$ rises from 0.  The interval of values of 
$\xi_j$ such that $x_i+\xi_j[\pi_j,-1]\in X_i$ lies in $[0,x_{ij}]$ (and
equals it unless $j$ is an essential good).  Hence by (2.7) there will be a 
unique maximizing quantity
$$
 \xi_j^\pls(x_i,\pi_j) = \argmax_{\xi_j}\! \Lset
      u_i(x_i +\xi_j[\pi_j,-1]\,)\Mset     
      x_i +\xi_j[\pi_j,-1]\in X_i\Rset. 
\eqno(2.14)$$
It will be either $x_{ij}$ or the unique value of $\xi_j$ where
$\theta'_\pls(\xi_j) =0$, whichever is reached first as $\xi_j$ increases,
and only starting then will selling cease to be attractive.  

The consequence for agent $i$ of selling that optimal amount would be to
replace $x_i$ by $x'_i= x_i +\xi_j[\pi_j,-1]$ for 
$\xi_j=\xi_j^\pls(x_i,\pi_j)$ and thereby change the price threshold 
$p_{ij}(x_i)$ for good $j$ to $p_{ij}(x'_i)$.  The case in (2.14) where 
optimality is attained with $\theta'_\pls(\xi_j)=0$ corresponds to having 
$p_{ij}(x'_i)=\pi_j$, in contrast to $p_{ij}(x_i)<\pi_j$.  The case where 
it is attained by hitting supply constraint for good $j$ while $\theta'_\pls$ 
is still increasing, corresponds instead to having $p_{ij}(x'_i)<\pi_j$.   
To summarize,
$$\eqalign{
\text{for the resultant holdings $x'_i$ at optimality in (2.14),} &\cr
\text{ $p_{ij}(x'_i) = \pi_j$ if $x'_{ij}>0$, but 
       $p_{ij}(x'_i) \leq \pi_j$ if $x'_{ij}=0$.}  &}
\eqno(2.15)$$
Note, however, that the case of $x'_{ij}=0$ can't come up if $j$ is an 
essential good for agent $i$, namely $j\in J_i$.  This is guaranteed 
by (2.7). 
\medskip

{\bf Agents as buyers.}  For an agent $i$ and a good $j\neq 0$, what needs 
to be considered in contemplating the purchase of a quantity $\xi_j>0$ 
of that good at a price $\pi_j>0$?  This consists of adding $\xi_j$
to the amount of good $j$ currently held while subtracting the money amount 
$\pi_j\xi_j$ from the current holding of money, namely $x_{i0}$, leaving
some remainder to 0.  It shifts $x_i$ to $x_i -\xi_j[\pi_j,-1]$, where the 
notation is again that of (2.9).  We are in a situation similar to that of 
a seller, but are looking this time at
$$
     \theta_\mns(\xi_j)= u_i(x_i -\xi_j[\pi_j,-1]\,)
\eqno(2.16)$$
for which
$$
  \theta'_\mns(\xi_j) = 
            -\frac{\partial u_i}{\partial x_{i0}}(x_i -\xi_j[\pi_j,-1]\,)\pi_j
            + \frac{\partial u_i}{\partial x_{ij}}(x_i -\xi_j[\pi_j,-1]\,),
\eqno(2.17)$$
For buying to be attractive, starting from an arbitrarily small amount
$\xi_j$, the marginal utility $\theta'_\mns(0)$ should be $>0$.  This time, 
however,
$$
   \theta'_\mns(0) = 
    \frac{\partial u_i}{\partial x_{i0}}(x_i)\Big[p_{ij}(x_i)-\pi_j \Big].
\eqno(2.18)$$ 
Therefore, 
$$\eqalign{
\text{buying good $j$ at price $\pi_j$ will be attractive to agent}&\cr 
\text{$i$ with current holdings $x_i$ if and only if 
     $\,\pi_j < p_{ij}(x_i)$.} &}
\eqno(2.19)$$

When buying is attractive in this way, how much should agent $i$ consider
buying?   Again the function $\theta_\mns(\xi_j)$ is concave with
$\theta'_\mns(\xi_j)$ strictly decreasing in $\xi_j$ and 
$\theta''_\mns(\xi_j)<0$, and the amount of money available for the
purchase is $x_{i0}>0$.  That amount can't end up being totally committed,
though, because we know from (2.7) that the interval of $\xi_j$ values
such that $x_i-\xi_j[\pi_j,-1]\in X_i$ and has utility at least that of
$x_i$ is a compact interval within $[0,x_{i0})$.  Thus, buying will continue 
to be attractive until reaching the level of
$$
 \xi_j^\mns(x_i,\pi_j) = \argmax_{\xi_j}\!\Lset
      u_i(x_i -\xi_j[\pi_j,-1]\,) \Mset x_i -\xi_j[\pi_j,-1]\in X_i\Rset,
\eqno(2.20)$$
which will be the unique value of $\xi_j$ yielding $\theta'_\mns(\xi_j)=0$.

The characterization of the optimality in (2.20) is more elementary than
it was for (2.14), where we had to cope with the possibility that the entire 
supply of good $j$ might be sold.  Here, if full purchase went through,, the 
new holdings vector would be $x'_i= x_i -\xi_j[\pi_j,-1]$ for 
$\xi_j=\xi_j^\mns(x_i,\pi_j)$, which comes unambiguously from achieving 
$\theta'_\mns(\xi_j)=0$.  The new price threshold $p_{ij}(x'_i)$ for good $j$ 
that replaces $p_{ij}(x_i)$ would emerge only out of that, hence
$$
\text{for the resultant holdings $x'_i$ at optimality in (2.20),}
           p_{ij}(x'_i) = \pi_j.
\eqno(2.21)$$

{\bf Utility function status.}  Like the thresholds $p_{ij}(x)$ in the 
conditions (2.13) and (2.19) that open the way to selling and buying, the 
quantities indicated in (2.14) and (2.20) for the amounts desired to be 
sold or bought {\it depend only on the preference relation, not the 
particular utility function $u_i$ with its scaling}.  That's because the 
optimization behind those quantities refers to proceeding along a line to a 
point of highest utility, and that point doesn't change with a change in 
scaling.  Agents could be tasked with locating that point ``geometrically'' 
with respect to the preference sets in the relation, but it's nicer to think 
of them working with a utility function $u_i$, as described, which entails 
no loss of generality.
\medskip

{\bf Connections with other work.}    Mandel and Gintis \cite{Mandel}
propose in their game model that agents have ``private'' vectors of prices
which can be deployed strategically, but they offer no basis in marginal 
utility, such as we provide here.  In the literature on strategic market 
games mentioned in our introduction, agents are expected to come up with 
bids for buying or selling up to some amount of a good, but the origin of 
such bids is left open.  Such price-quantity pairs can now be seen as 
coming, for instance, from our analysis of potential buyers and sellers in 
association with their price thresholds.

\section{ 
          Bilateral trading in one good at a time}

With the characteristics and interests of sellers and buyers in hand, we
can proceed with ideas about how a seller and a buyer might come together 
through information revealed by one or the other, or both, and proceed
with a money-based transaction in a good $j\neq 0$.
\medskip

{\bf Price premiums and bid-ask spreads.}
It has been seen that an agent $i$ with threshold price $p_{ij}(x_i)$ for 
good $j$ only wants to sell at a higher price than that, and only wants to 
buy at a lower price.  How should those prices be set?  Our approach is 
that agent $i$ has at any given stage a ``premium'' $\delta_{ij}>0$ in mind 
to make a transaction in good $j$ worth undertaking.  The selling and buying prices are set by the agent to maintain that premium:
$$
     p_{ij}^\pls(x_i) = p_{ij}(x_i) +\delta_{ij},\qquad 
     p_{ij}^\mns(x_i) = p_{ij}(x_i) -\delta_{ij}.
\eqno(3.1)$$
The particular premiums are only temporary.  Eventually they will need to 
be reduced more and more if trading is going to be able to identify an 
equilibrium to utmost precision.

Agent $i$ can reveal being open to selling a good $j$ at price 
$p_{ij}^\pls(x_i)$ or being open to buying it at price $p_{ij}^\mns(x_i)$, 
or both, or neither --- in the context of holding off until later.  Just 
how this may be articulated over time will be the subject of discussion 
later, in the Section 4.

\proclaim Trading action based on proposed prices.   A bilateral trade is 
available in good $j$ between agent $i_1$ as seller and agent $i_2$ as 
buyer at a price $\pi_j$ if, with respect to (3.1),
$$  
   x_{i_1 j}>0, \qquad  p_{i_1 j}^\pls(x_{i_1}) \leq p_{i_2 j}^\mns(x_{i_2}), 
     \qquad \pi_k \in [ p_{i_1 j}^\pls(x_{i_1}),p_{i_2 j}^\mns(x_{i_2})] 
\eqno(3.2)$$
In acting on this availability (if they choose to do so), the amount of
good $j$ transferred from the seller to the buyer will be
$$
   \xi_j = \min\{ \xi_j^\pls(x_{i_1},\pi_j),\,\xi_j^\mns(x_{i_2},\pi_j)\} 
\eqno(3.3)$$
as determined from (2.14) and (2.20).  In return for this amount of good
$j$, the buyer will pay the seller the money amount $\pi_j\xi_j$.

It would not really be necessary to trade all the way up to the 
profitability-limiting level in (3.3), but this rule is the simplest.
It guarantees that the agents will have activated at least one of the 
optimality rules in (2.14) or (2.20), and that will be useful to us later. 

Along with these specifics about how agents can mutually improve the
utility values of their holdings by trading individual goods, it will be
important to understand the situation in which the limits of such trading 
have been reached.   It turns out that this is the case where an 
equilibrium of prices and holdings is at hand.

\proclaim Theorem 1 {\rm (seller-buyer characterization of an equilibrium
of prices and holdings)}.
Suppose the agents have holdings $\bar x_i$ such that no bilateral trade
is available in any good at any premium levels $\delta_{ij}>0$,
or in other words:  
$$\eqalign{
\text{for every good $j\neq 0$ and potential seller $i_1$ 
                 and buyer $i_2$ } &\cr
\text{enabled by a supply  $\bar x_{i_1j}>0$, one has} 
          p_{i_1j}(\bar x_{i_1})\geq p_{i_2j}(\bar x_{i_2}).  &}
\eqno(3.4)$$
Then the price vector 
$$
   \bar p = (\bar p_1,\ldots,\bar p_n) \text{with}
   \bar p_j = \max_{i=1,\ldots,m} p_{ij}(\bar x_i)>0
\eqno(3.5)$$
for the nonmonetary goods (the money price of money being 1) furnishes an 
equilibrium with those holdings: 
$$\eqalign{
 \;\, \bar x_i \text{maximizes} \,u_i(x_i)\, 
    \text{over} x_i\in X_i \text{subject to the} &\cr
   \text{budget constraint} x_{i0}+\sumn_{j=1}^n \bar p_j x_{ij} 
         = \bar x_{i0}+\sumn_{j=1}^n \bar p_j \bar x_{ij}. &}
\eqno(3.6)$$
Conversely, every equilibrium of prices and holdings can be described in 
this manner.

\state Proof.  If the threshold inequality at the end of (3.4) were strict, 
and $\bar x_{i_2 j}>0$ as well, one could consider reversing the roles of 
buyer and seller and thereby opening up a trade.  Thus (3.4) really says 
for each good $j\neq 0$ that
$$\eqalign{
   \text{$p_{ij}$ has the same value for all agents $i$ 
                     having $\bar x_{ij}>0$,} &\cr
   \text{whereas $p_{ij}$ is $\,\leq\,$ that value for agents $i$ 
                     with $\bar x_{ij}=0$.} &}
\eqno(3.7)$$
Storing this observation temporarily, let's next analyze the optimality 
in (3.6) from the perspective of convex analysis.  In terms of the
linear function
$$
      l_i(x_i)= x_{i0}+\sumn_{j=1}^n \bar p_j x_{ij} =x\mdot (1,p)
\eqno(3.8)$$
(3.6) says $\bar x_i$ maximizes $u_i(x_i)$ over $\reals^n_\pls$ 
subject to $l_i(x_i)\leq l_i(\bar x_i)$, but that's the same as saying 
$\bar x_i$ minimizes $l_i(x_i)$ over the set 
$C_i=\lset x_i\geq 0 \mset u_i(x_i)\geq u_i(\bar x_i)\rset$.  
That set is convex, so this holds if and only if the normal cone 
$N_{C_i}(\bar x_i)$ to $C_i$ at $\bar x_i$ contains the vector 
$-\nabla l_i(\bar x_i)=-(1,p)$.  The normal cone is the sum of the ray 
generated by $\nabla u_i(\bar x_i)$ and the normal cone to $\reals^n_\pls$
at $\bar x_i$, consisting of all $v\leq 0$ such that $v\mdot\bar x_i=0$.  
Having $-(1,p)$ belong to it corresponds therefore to the existence of 
$\lambda_i> 0$ such that
$$
    \nabla u_i(\bar x_i)-\lambda_i(1,p)\geq 0
   \text{with equality in the $j$\/th component if $\bar x_{ij}>0$, }
\eqno(3.9)$$
where in particular $\bar x_{i0}>0$ as enforced by (2.4)--(2.5).
Component by component, this condition that is equivalent to (3.6) means 
$$
   \frac{\partial u_i}{\partial x_{ij}}(\bar x_i)\left\{\eqalign{ 
        =\lambda_i \text{for $j=0$,} &\cr
        =\lambda_i\bar p_j \text{for $j\neq 0$ with $\bar x_{ij}>0$,}&\cr
     \leq\lambda_i\bar p_j \text{for $j\neq 0$ with $\bar x_{ij}=0$.}&}
     \right. 
\eqno(3.10)$$
Another way of stating it, in terms of the definition of price
thresholds, and understanding that the first line of (3.12) identifies 
$\lambda_i$ as the money partial derivative 
$\frac{\partial u_i}{\partial x_{i0}}(\bar x_i)$, is that
$$
    p_{ij}(\bar x_i)\left\{\eqalign{ 
        =\bar p_j \text{for $j\neq 0$ with $\bar x_{ij}>0$,}&\cr
     \leq\bar p_j \text{for $j\neq 0$ with $\bar x_{ij}=0$.}&} \right. 
\eqno(3.11)$$
But this version of (3.6) obviously is the combination of (3.7) and (3.5),
and we are done.  \eop
     
The new and valuable conclusion here is that {\it bilateral trades 
involving bundles of several goods at a time aren't required in support 
of equilibrium}, not to speak of multilateral trades in which several agents 
are simultaneously engaged.  Feldman already in 1973 in \cite{Feldman} 
showed that pairwise Pareto optimality, in which no agent pair can jointly 
improve utility by a mutual exchange of goods, implies full Pareto 
optimality if some good is always attractive to every agent --- like money.
But this was with multigood trades. 

Trading as we depict it conserves goods;  after each transaction the total 
of each good $j=0,1,\ldots,n$ held by the agents agrees with the amount 
that was present initially.  Thus, the equilibria potentially reachable 
from total initial supplies $s_j$ correspond to the holdings vectors 
$\bar x_i$ and prices $\bar p_j$ satisfying (3.11) along with 
$$
         \sumn_{i=1}^m \bar x_{ij}=s_j \text{for} j=0,1,\ldots,n.
\eqno(3.12)$$
These conditions describe the ``equilibrium manifold'' associated with the
given utility functions $u_i$ and supply vector $s=(s_0,s_1,\ldots,s_n)
\in \reals_{\pls\pls}^n$ with respect to the $m(n+1)+n$ unknowns
$\bar x_{ij}$ and $\bar p_j$.  When all goods are essential, the inequality 
case drops out of (3.11), leaving $mn$ equations in the unknowns to be
combined with the $n+1$ equations in (3.12), so that the number of
equations exceeds the number of unknowns by $m-1$.  If degeneracy doesn't 
intervene in the equations, they ought to produce a kind of 
$(m-1)$-dimensional ``surface'' as the manifold in question.  {\it This 
nonuniqueness is what we expect and embrace}.  We aren't in the Walras 
framework and are instead exploring an approach to equilibrium based 
squarely on buying and selling as ordinarily envisioned, instead of relying 
on an abstract societal entity to coordinate the agents' wishes and engineer 
a simultaneous grand exchange of goods among them.

\section{
             Evolving toward equilibrium}

The central question we now wish to answer, at least to some level of
satisfaction, is whether, under additional specifics and assumptions if 
necessary, the kind of direct trading that has been described will lead
in the limit to equilibrium holdings as characterized in Theorem~1.

Getting back first to pinning down how a seller and buyer may 
come together and come up with a price $\pi_j$ for a good $j\neq 0$, we 
can foresee several possibilities.  After a seller $i_1$ has proposed a 
price $p_{i_1 j}^\pls$, a buyer $i_2$ may respond by simply accepting that 
price as $\pi_j$, as long as $p_{i_1 j}^\pls\leq p_{i_2 j}^\mns$.  In that 
case the buyer acts without revealing $p_{i_2 j}^\mns$.  Likewise, a 
seller $i_1$ could respond to a price $p_{i_2,j}^\mns$ proposed by a 
buyer $i_2$ and accept that as $\pi_j$, as long as 
$p_{i_1 j}^\pls\leq p_{i_2, }^\mns$, and thereby act without revealing 
$p_{i_1 j}^\mns$.   More complicatedly, if both $i_1$ and $i_2$ have 
revealed prices acceptable to them, and those prices are different, 
$p_{i_1 j}^\pls < p_{i_2 j}^\mns$, they could take $\pi_j$ to be the 
average, say.  Many details could obviously be added also to explain how 
agents find each other.  There might be ``information boards'' or even 
trading posts with double auctions in various goods --- whatever could
work in the end to effect elementary exchanges.  

The important thing, regardless, is that each time a bilateral trade occurs 
in some good the utility values of the holdings of both agents in the trade
are improved (while those of the other agents are unaffected).  We envision 
a scheme in which prices based on premiums are proposed by buyers and 
sellers to an extent that, in successive iterations, all agents and goods 
are repeatedly brought into the picture.  Because utility values of the 
holdings of the agents go up each time for seller and buyer, and none ever 
go down, it's impossible that the configuration of holdings can ``recycle.''  
After each iteration it will be different, and in a Pareto sense always 
better.

Our hope is that eventually in this manner a stage will be reached in 
which no bilateral trade in any good is available --- at the current 
premium levels $\delta_{ij}$.   At that point the premiums can be lowered 
and the process restarted.\footnote{
    Of course this could be triggered good-by-good instead.}  
But how do we know that the iterations won't stagnate through utility
improvements falling short of what they ought to be?   

To head off such stagnation, we'll introduce a restriction on the scope 
of our goods model in (2.3), namely that
$$
   \text{henceforth {\it all\/} goods are essential, so that
                        $X_i=\reals^n_{\pls\pls}$ for all agents $i$.}
\eqno(4.1)$$
Then, as explained at the end of Section 3, we anticipate facing an
equilibrium manifold of dimension $m-1$ as the target to be reached 
somewhere by the trading scheme.  Of course, the particular equilibrium 
that might take shape from the bilateral activity will likely be affected 
by the random order in which agents come together and choose goods to 
look at.

The key consequence of assumption (4.1), which is common to much of the
literature on equilibrium, will be to provide quantitative bounds on the 
sizes of the improvements in utility that result from trading.  Something 
weaker than (4.1) might well be possible through further research.  An easy 
step toward relaxation will be explained at the end of this section, in
particular.

\proclaim Theorem 2 {\rm (lower bounds on improvements from trading)}.
Under the utility assumptions strengthened by (4.1) there exists $\mu>0$, 
independent of anything other than the utility functions and the agents'
initial holdings, for which the following property holds.  When a trade in 
good $j$ takes place between a seller $i_1$ with premium level 
$\delta_{i_1 j}$ and a buyer $i_2$ with premium level $\delta_{i_2 j}$, 
their utility levels will improve to at least the degree that, when the 
current holdings $x_{i_1}$ and $x_{i_2}$ are replaced by the after-trade 
holdings $x'_{i_1}$ and $x'_{i_2}$, one will have
$$
  u_{i_1}(x'_{i_1})\geq u_{i_1}(x_{i_1})+\mu\delta_{i_1 j}^2, \qquad
  u_{i_2}(x'_{i_2})\geq u_{i_2}(x_{i_2})+\mu\delta_{i_2 j}^2.
\eqno(4.2)$$

\state Proof.  Let $s_j>0$ denote the total supply of good $j$, which 
starts as $\sumn_{i=1}^m x_{ij}^0$ and is forever maintained under trading.
Let $s=(s_0,s_1,\ldots,s_n)$.  The holdings $x_i$ of agent $i$ will always 
then lie in the compact convex set
$$
      C_i = \lset x_i \mset u_i(x_i)\geq u_i(x_i^0),\; x_i \leq s \rset
     \subset \reals^{n+1}_{\pls\pls}.
\eqno(4.3)$$
Because $u_i$ is a twice-continuously differentiable function on 
$\reals^{n+1}_{\pls\pls}$, its first and second partial derivatives are 
bounded on $C_i$.   In particular, by virtue of lower and upper bounds on 
the first partial derivatives, these necessarily being positive, there will 
exist positive lower and upper bounds on the price thresholds $p_{ij}(x_i)$ 
in (2.4) as $x_i$ ranges over $C_i$.  Since any price $\pi_j$ that might
enter a transaction for good $j$ has to lie between a seller's threshold
and a buyer's threshold, such prices are thus limited to some interval
$[\pi_j^\mns,\pi_j^\pls]\subset (0,\infty)$.

Consider now the function $\theta_\pls(\xi_j)$ in (2.10) with notation 
$[\pi_j,-1]$ as explained in (2.9).   This is the utility expression that 
agent $i$ as a seller would be involved with optimizing over a line 
segment necessarily lying within $C_i$.  We have already calculated its
derivative at $0$ in (2.12) as 
$\frac{\partial u_i}{\partial x_{i0}}(x_i)[\pi_j-p_{ij}(x_i)]$, but here
$\pi_j-p_{ij}(x_i) \geq\delta_{ij}$, so can say that
$$
    \theta'_\pls(0) \geq \alpha_i \delta_{ij}
    \text{where $\alpha_i>0$ is a lower bound for} 
     \frac{\partial u_i}{\partial x_{i0}}(x_i) \text{on $C_i$.}
\eqno(4.4)$$
The second derivative of $\theta_\pls(\xi_j)$ can be calculated in turn 
from (2.10) as
$$
   \theta''_\pls(\xi_j)= [\pi_j,-1]\mdot\nabla^2 u_i(x'_i)[\pi_j,-1]
                 \text{for} x'_i=x_i+\xi_j[\pi_j,-1].
\eqno(4.5)$$
Because $\pi_j$ must be in the interval $[\pi_j^\mns,\pi_j^\pls]$, while 
$x_i$ and $x'_i$ must lie in $C_i$, compactness and continuity ensure that 
the values in (4.5) over all such possible arguments are bounded from below.
Thus, there exists $\beta_{ij}>0$ such that 
$$
     \theta''_\pls(\xi_j)\geq -\beta_{ij} \text{for all feasible $\xi_j$,}
\eqno(4.6)$$
where the feasibility refers to $\xi_j$ being admissible in any transaction
that agent $i$ might enter in selling good $j$.  The combination of (4.4) 
and (4.6) implies that
$$
 \theta'_\pls(\xi_j)\geq \alpha_i\delta_{ij}-\beta_{ij}\xi_j, 
  \qquad
 \theta_\pls(\xi_j) \geq \theta_\pls(0) +\alpha_i\delta_{ij}\xi_j
                            -\half \beta_{ij}\xi_j^2,
\eqno(4.7)$$
for all feasible $\xi_j$.  Therefore
$$
   \max_\substack{\xi_j\geq 0}{x_i+\xi_j[\pi_j,-1] \in X_i} \theta_\pls(\xi_j) 
   \geq \max_\substack{\xi_j\geq 0}{x_i+\xi_j[\pi_j,-1] \in X_i}
  \lset \theta_\pls(0)+ \alpha_i\delta_{ij}\xi_j -\half\beta_{ij}\xi_j^2\rset, 
\eqno(4.8)$$
as well as
$$
    \argmax_\substack{\xi_j\geq 0}{x_i+\xi_j[\pi_j,-1] \in X_i} 
       \theta_\pls(\xi_j) \geq
    \argmax_\substack{\xi_j\geq 0}{x_i+\xi_j[\pi_j,-1] \in X_i}\lset 
   \theta_\pls(0)+ \alpha_i\delta_{ij}\xi_j-\half\beta_{ij}\xi_j^2\rset, 
\eqno(4.9)$$
where, due to (4.3), each maximum is obtained with a value of $\xi_j$ such 
that $x_i+\xi_j[\pi_j,-1]$ belongs to the interior of the goods orthant, 
hence simply by setting the derivative equal to 0.  On the 
right of (4.9), that leads to specifically to
$$
   \xi_j =\frac{\alpha_i}{\beta_{ij}}\delta_{ij}, \text{yielding max value}
   \theta_\pls(0)+\frac{\alpha_i}{2\beta_{ij}}\delta_{ij}^2.
$$
The conclusion reached is that the new holdings vector $x'_i$ obtained by
agent $i$ as seller will have 
$$
     u_i(x'_i)\geq u_i(x_i) +\mu_{ij}\delta_{ij}^2 \text{for}
            \mu_{ij} = \frac{\alpha_i}{2\beta_{ij}}>0.
$$

Treating agent $i$ as buyer instead of seller leads to the same end
through analysis of the function $\theta_\mns(\xi_j)$ in (2.12), because
all that really changes is that $[\pi_j,-1]$ is replaced by $-[\pi_j,-1]$, 
and that has no effect on the hessian expression for second derivatives.  
Taking $\mu$ to be the smallest of the factors $\mu_{ij}$ over all $i$ and
$j$ justifies the double improvement claim in (4.2).  \eop

A definitive answer to the question of convergence to an equilibrium is
now at hand for our scheme in which agents interact on the basis of 
particular price premiums $\delta_{ij}>0$ as long as trades are possible,
but reduce those premiums once trades are blocked, and keep doing that as
the premium levels tend to 0.

\proclaim Theorem 3 {\rm (convergence to an equilibrium)}.
Under the additional assumption in (4.1), the outlined trading scheme for 
agents at a given level of price premiums $\delta_{ij}>0$ can proceed for 
only finitely many trades before no further trades are available.  If the 
trading is restarted then at lower levels of premiums, and the 
process repeats iteratively with still lower levels, tending to 0, the 
holdings $x_i$ and price thresholds $p_{ij}(x_i)$ converge to some 
particular equilibrium of prices and holdings as characterized in 
Theorem 1, moreover with
$$
   \text{all goods $j\neq 0$ having} p_{ij}(\bar x_i) = \bar p_j 
     \text{for every agent $i$.}
$$

\state Proof.  The trading scheme, however articulated among the traders, 
generates for each agent $i$ a sequence of holdings vectors $x_i$ within 
the compact set $C_i$ in (4.3).  For a cluster point $\bar x_i$ of that 
sequence, $p_{ij}(\bar x_i)$ will be a cluster point of the corresponding 
price thresholds $p_{ij}(x_j)$ for each good $j$.  We claim two things.  
First, $p_{ij}(\bar x_j)$ will have the same value for every agent $i$, so 
that this value, as $\bar p_j$ will provide prices needed along with the 
holdings $\bar x_i$ to constitute an equilibrium as described in Theorem 2.  
And second, that the cluster points are unique:  even though the existence 
of more than one equilibrium isn't excluded, the sequences that are 
generated converge to particular limits.
 
The process under fixed premiums results in improvements as in Theorem 2 
after each transaction.  There is a lower bound that way to the size of 
the utility improvements.  Because the utility functions $u_i$ are bounded 
from above over the set $C_i$, in which holdings always remain, stagnation 
is impossible.

Once the stage is reached when no further transactions are possible at the
current premium levels $\delta_{ij}$, we must have for every potential
seller $i_1$ and buyer $i_2$ that $p_{i_1 j}^\pls > p_{i_2,j}^\mns$ for
the acceptable prices in (3.1).   This means that
$
     p_{i_1 j}(x_{i_1})+\delta_{i_1 j} > p_{i_2 j}(x_{i_2})-\delta_{i_2 j},
$
but since all agents have positive holdings in all goods in consequence of
(4.7), the roles of seller and buyer can be interchanged.  Thus, we must have
$$
  |p_{i_1 j}(x_{i_1})- p_{i_2 j}(x_{i_2})| < \delta_{i_1 j} +\delta_{i_2 j}
       \text{for all pairs of agents $i_1$ and $i_2$.}
\eqno(4.10)$$
Clearly then, as the trading restarts and proceeds over and over with
lower premium levels tending to $0$, it must be that, for the subsequence
of holdings $x_i$ that converges to $\bar x_i$ at the beginning of this 
proof, we end up for each agent $i$ and good $j$ with the same value 
for $p_{ij}(\bar x_i)$.

For the second claim, we recall that levels of utility only increase for
the agents as  the trading proceeds, never decrease.  The vector of those 
levels thus tends upward to a vector $(\bar u_1,\ldots,\bar u_m)$.  Any 
array $[\bar x_1,\ldots,\bar x_m]$ of equilibrium holdings that is
approached as a cluster point of the (bounded) sequence of holdings 
vectors generated by the trading must have 
$$
      u_i(\bar x_i) = \bar u_i  \text{for} i=1,\ldots,m. 
\eqno(4.11)$$
Suppose that, along with  $[\bar x_1,\ldots,\bar x_m]$, there is another 
array $[\bar x'_1,\ldots,\bar x'_m]$ that is an equilibrium cluster point,
likewise then having 
$$
      u_i(\bar x'_i) = \bar u_i  \text{for} i=1,\ldots,m. 
\eqno(4.12)$$
Because total supplies of goods are maintained in the process, we know
that
$$
    \sumn_{i=1}^m \bar x_i = \sumn_{i=1}^m \bar x'_i= \sumn_{i=1}^m x_i^0.
\eqno(4.13)$$
Let $x_i^*=\half[\bar x_i+\bar x'_i]$, so that again
$$
    \sumn_{i=1}^m x_i^* = \sumn_{i=1}^m x_i^0.
$$
Because each $u_i$ is strictly (even strongly) concave along line segment
joining $\bar x_i$ and $\bar x'_i$ inasmuch as the difference is a mixed-sign
vector when nonzero by (4.13), and some differences are indeed nonzero,
we have
$$
   u_i(x_i^*) \geq \half u_i(\bar x_i)+\half u_i(\bar x'_i) =\bar u_i
       \text{for all $i$, with strict inequality for some $i$.}
$$ 
In view of (4.11) and (4.12), this says that neither
$[\bar x_1,\ldots,\bar x_m]$ nor $[\bar x'_1,\ldots,\bar x'_m]$ is Pareto
optimal.  But Pareto optimality holds for any exchange equilibrium as a
special case of it holding for any Walras equilibrium.  \eop   

The achievement of an equilibrium as a limit and not just as a cluster
point is a unique contribution here.   In other work, such as that of
Flaam in \cite{Flaam1}, only a cluster point level of convergence is
obtained, if anything is obtained at all.   
\medskip

{\bf Remark about premiums.}  
The scheme in Theorem 3 of having the agents progressively 
lower their price premiums when trading opportunities are lacking is
important for the purpose of proving eventual convergence.  From a practical
perspective in economics, however, it might be more natural to imagine 
these money-denominated premiums as characteristics of the agents along 
with their utility functions.  With their levels fixed, buying and selling 
would come to a halt short of a full equilibrium.  But this would reflect 
reality.   Money prices anyway are always truncated in practice.  The 
smallest money unit ever contemplated in market transactions could serve 
for instance uniformly as every premium $\delta_{ij}$.
\medskip

{\bf Concluding comments about the utility assumptions}.
The convergence in Theorem 3 rests on the lower bounds for improvements
that were derived in Theorem 2, and this is where our second-order
assumptions on the utility function $u_i$ and its concavity are most
important.  How plausible are those assumptions from the angle of
economics?   In assessing this, we have background support from 
\cite[Theorem 4.2]{Preferences} in regarding the representability of a 
preference relation by a $\cC^2$ quasi-concave utility as a property of the
basic workings of that the relation, not a wishful construct.  It can well
be a natural feature of the agents' preferences in our model.  Local strong 
convexity of the preference sets can be interpreted like that as well.

It helps now to keep in mind that the agents in our setting will not make 
use of the entire goods orthant but only a limited portion indicated by
the compact sets $C_i$ in (4.3).  According to 
\cite[Theorem 3.5]{Preferences}, a strongly convex preference relation for 
agent $i$ that's open to representation by a $\cC^2$ utility function $u_i$ 
has one which is concave relative to the set $C_i$, and minimally so, 
with it then being unique up to just the choice of units in which utility 
is to be measured.  

That puts everything on solid ground.  But it might still be wondered 
whether an agent could know such a function and be able to deploy it --- 
pure existence might not be enough.   In our scheme, though, such knowledge 
actually isn't required!  We pointed out at the end of Section 2 that 
the optimization steps in selling and buying can, in principle, be carried 
out directly in terms of the preference relation.  The utility 
representation essentially serves only for convenience.  And as for the 
strong concavity yielding the bounds in Theorem 2, it merely serves in the 
proof of convergence through its existence.  An agent doesn't have to do
anything to take advantage of it or even be aware of it.
\medskip

\goodbreak
{\bf An easy extension.}
The extra assumption in (4.1) that we made in order to get the improvement
guarantee in Theorem 2 as a key ingredient of the convergence proof in
Theorem~3 makes every agent to insist on holding always a positive
amount of every good.  That's an unpleasant requirement, especially in
view of the preceding analysis of buying and selling, which didn't need it.
It's good to note, therefore, that (4.1) can be replaced by a somewhat
weaker assumption which permits agents to hold zero quantities of goods in
which they have no interest in at all.  For that, the sets $X_i$ in (2.3) 
can be modified from the beginning to 
$$\eqalign{
  X_i= \lset x_i\in\reals^{n+1}_\pls \mset 
     x_{ij}>0 \,\txt{for} j\in J_i,\; 
     x_{ij}=0 \,\txt{for} j\in J^0_i \rset &\cr \hskip35pt  
  \text{for disjoint collections $J_i$ and $J_i^0$ of goods $j$.} &}
\eqno(4.13)$$
The utility assumptions can be altered so as not to apply to the goods
$j\in J_i^0$, which have no effect on them.  Everything continues to work
in this mode with only minor and obvious adjustments to the story.

\section{ 
              Computational experiments}

The effectiveness of the bilateral trading dynamic in Sections 3 and 4 is 
open to being explored by computer simulation.  Examples of that will now 
be presented.  After giving a formal statement of the trading procedure as
Algorithm 1, we apply it first to situations in which agents have utility 
functions in the Cobb-Douglas family and the holdings of all goods are 
necessarily positive, in accordance with assumption (4.1).  The conditions 
behind Theorem 3 are met, and numerical convergence to a market equilibrium 
is consistently obtained, as expected.  

\begin{algorithm}
	\caption{Bilateral Trading Scheme}
	\begin{algorithmic}[1]
		\Procedure{BTS}{Agent parameters 
    $\{\beta_{ij}\},\,\{x_{ij}^0\},\{\delta_{ij}^0\}$, Algorithm parameters 
    $\epsilon_\delta,\,\epsilon_p,\lambda$}   
		\State Define $p^0_{ij}=p_{ij}(x_{i}^0),\forall i,j$.
		\For{$k=1,\ldots$}
		\If{$\max_j\operatorname{stdev}{(p_{ij}^k)}<\epsilon_p$} 
		\Return {\bf Equilibrium found}.
		\ElsIf{$\max_{ij}\{\delta^k_{ij}\}<\epsilon_\delta$} 
		\Return Max of deltas is less than the given tolerance 
                                 level of $\epsilon_\delta$
		\Else
		\For{$i_1\in\Pi(1,\ldots,m)$} 
		\For{$i_2\in\Pi(1,\ldots,m)$} 
		\For{$j\in\Pi(1,\ldots,n)$}
		\If{Trade between agents agents $i_1$ and $i_2$ of good $j$ 
                                        is possible}
		\State Update holdings $x_{i_1,j}^{k+1}$, 
                                                   $x_{i_2,j}^{k+1}$
		\State Update agents' prices $p^{k+1}_{i_1}$, 
                                                 $p^{k+1}_{i_2}$.
		\State Go to Step 4  
		\EndIf
		\EndFor
		\EndFor
		\EndFor
		\If{No trade is possible}
		\State $\delta^{k+1}_{ij}\leftarrow \lambda \delta^k_{ij}$, 
                  make $k \leftarrow k+1$ and go to Step 7
		\EndIf
		\EndIf
		\EndFor
		\EndProcedure
	\end{algorithmic}
\end{algorithm}

We then experiment with whether the trading procedure can reach an 
equilibrium even in a example where the interiority assumption (4.1) is 
not in force, so that agents can hold zero quantities of some goods.  In 
that example the utility functions are not Cobb-Douglas.  Convergence to
an equilibrium is sometimes obtained but can also bog down.  The lesson
is that a refined substitute for the utility growth guarantee in
Theorem 2 is definitely needed for convergence outside of assumption (4.1) 
or its relaxation described at the end of the preceding section. 

In the inspection process for determining the existence of bilateral 
trade opportunities, we rely here on a naturally {\it random trading order}.
This means that we inspect potential buyer/seller partners at random, 
traveling through all possible combinations and for each one looking 
randomly through all the goods to possibly trade, but stopping at the first 
trading opportunity encountered.  If the inspections comes to an end with 
no trading opportunity at all having emerged, the agents all reduce their 
premium levels $\delta_{ij}$, and at those new levels the inspection 
resumes.  The ultimate convergence of this scheme is guaranteed by Theorem 3.  
In Steps 7, 8 and 9 of Algorithm 1 the functions $\Pi$ represent the 
generated permutations of agents or goods.

When the procedure terminates in Step 4 in iteration $k$ with an equilibrium 
having been found --- approximately --- under the stopping criterion, the
prices and holdings associated with it are taken to be
$$
   \bar p_j =\frac{1}{m}\sumn_{i=1}^m p^k_{ij}, \qquad
       \bar x_{ij} = x^k_{ij}.
\eqno(5.1)$$

We developed a Python 3.8.11 code for this algorithm, making it freely 
available in our repository\footnote{
     Github repository  \href{https://github.com/jderide/DeJoRo-Bilateral}
          {https://github.com/jderide/DeJoRo-Bilateral} } 
along with the parameters and documentation.  We ran our examples on an 
Apple Mac Mini 2018 (Intel Core i5 CPU@3GHz, 32 GB memory), operating with
macOS Big Sur.

\proclaim Example 1.
This concerns an economy with $n=9$ non-money goods and $m=5$ agents, in 
which each agent $i$ has a utility function of Cobb-Douglas type: 
$$
       u_i(x_i)=\prod_{j=0}^n x_{ij}^{\beta_{ij}}, \text{with}
     0<\beta_{ij}<1, \quad\sum_{j=0}^n\beta_{ij}<1.
\eqno(5.2)$$
The values assigned to the exponents $\beta_{ij}$ are shown in Table 1 
along with the quantities taken as the agents' initial holdings 
$x_{ij}^0$, all of which are positive.  
	
	\begin{table}[h]  
		\begin{footnotesize}
			\centering
			\begin{tabular}{cc|cccccccccc}
				Agent& Parameter   & $j=0$ & $j=1$ & $j=2$ & $j=3$ & $j=4$ & $j=5$ & $j=6$ & $j=7$ & $j=8$ & $j=9$ \\
				\hline
				\multirow{2}{*}{$i=1$} & $\beta_{1}$ & 0.09  & 0.09  & 0.09  & 0.09  & 0.09  & 0.09  & 0.09  & 0.09  & 0.09  & 0.09  \\
				& $x_1^0$     & 59    & 76    & 10    & 37    & 54    & 99    & 73    & 30    & 25    & 20    \\
				\hline
				\multirow{2}{*}{$i=2$} & $\beta_{2}$ & 0.05  & 0.1   & 0.17  & 0.02  & 0.16  & 0.1   & 0.16  & 0.07  & 0.03  & 0.04  \\
				& $x_2^0$     & 14    & 40    & 63    & 57    & 69    & 39    & 34    & 86    & 10    & 56    \\
				\hline
				\multirow{2}{*}{$i=3$} & $\beta_{3}$ & 0.06  & 0.05  & 0.09  & 0.15  & 0.07  & 0.08  & 0.14  & 0.02  & 0.11  & 0.13  \\
				& $x_3^0$     & 19    & 57    & 43    & 65    & 78    & 40    & 9     & 82    & 71    & 82    \\
				\hline
				\multirow{2}{*}{$i=4$} & $\beta_{4}$ & 0.01  & 0.15  & 0.01  & 0.11  & 0.11  & 0.16  & 0.03  & 0.14  & 0.09  & 0.09  \\
				& $x_4^0$     & 10    & 65    & 35    & 43    & 63    & 74    & 79    & 38    & 20    & 27    \\
				\hline
				\multirow{2}{*}{$i=5$} & $\beta_{5}$ & 0.03  & 0.13  & 0.05  & 0.16  & 0.16  & 0.07  & 0.08  & 0.1   & 0.08  & 0.04  \\
				& $x_5^0$     & 37    & 70    & 40    & 94    & 83    & 15    & 34    & 97    & 35    & 34   
			\end{tabular}
	\caption{Utility parameters $\beta_{ij}$ and initial holdings 
                               $x^0_{ij}$ in Example 1}
		\end{footnotesize}
	\end{table}
	
We ran our algorithm 50 times with a random inspection strategy and observed 
that the trading scheme reached a market equilibrium every time in the
sense of the stopping criterion.  The precision level for that criterion
in Step 4 of the algorithm was taken to to have $\epsilon_p=10^{-6}$,, so
the equilibrium price vector $\bar p$ reported as indicated in (5.1) was 
such that its distance from every one of the agent price vectors $p_i^k$ was
less than $10^{-6}$.  The median execution time for these runs was 2.12 
seconds in 3093 iterations.  

Because of the random inspection strategy, each run could, and did, end
with a different equilibrium, and significant differences in the
resulting prices were observed.  Rather that reporting the results of all
50 runs, we selected the four with the highest absolute price differences,
namely runs $5,\,4,\,9,\,1$, for display of their equilibrium prices in 
Table 2.  

	\begin{table}[h]   
		\centering
       \begin{tabular}{l|rrrrrrrrrr}
	 &  $j=0$ &  $j=1$ &  $j=2$ &  $j=3$ &  $j=4$ &  $j=5$ &  $j=6$ &  $j=7$ &  $j=8$ &  $j=9$ \\
	\hline
	$\overline p^5$   & 1.0000 & 0.9278 & 1.2147 & 1.0510 & 0.9491 & 1.0213 & 1.2575 & 0.6828 & 1.4315 & 1.0215 \\
	$\overline p^4$   & 1.0000 & 0.9105 & 1.2176 & 1.0077 & 0.9356 & 0.9999 & 1.2394 & 0.6736 & 1.3860 & 0.9888 \\
	$\overline p^9$   & 1.0000 & 0.9468 & 1.2214 & 1.0493 & 0.9678 & 1.0303 & 1.2552 & 0.7005 & 1.4263 & 1.0077 \\
	$\overline p^1$   & 1.0000 & 0.9192 & 1.2000 & 1.0141 & 0.9387 & 1.0009 & 1.2229 & 0.6841 & 1.3907 & 0.9836 \\
	$\overline p^{W}$ & 1.0000 & 0.9575 & 1.2218 & 1.0569 & 0.9680 & 1.0594 & 1.2609 & 0.7102 & 1.4501 & 1.0371 \\
     \end{tabular}
	\caption{ Selected outcomes for equilibrium prices in Example 1}
	\end{table}	

The corresponding money-denominated prices $p^W_j$ for a Walras equilibrium 
with these utility functions and initial holdings are also shown in Table 2 
for comparison and possible side interest in underscoring a difference in aims 
and concept.  The Walrasian model {\it lacks true market legitimacy}, as 
explained in Section 1 and \cite{WalrasExchanges}.  There is no reason why 
a market equilibrium achieved through our bilateral trading scheme would 
reproduce its configuration of terminal holdings and prices, which are not
known to be supported by market activities in any ordinary sense, but
depend rather on the intervention of a clearing house agent.  

To get the prices $p^W_j$ in Table 2, we proceeded as follows.  In the 
Cobb-Douglas setting of Example 1 with its coefficients $\beta_{ij}$ 
and initial holding $x^0_{ij}$, the final holdings $\tilde x_{ij}$ and 
relative prices $\tilde p_j$ in the Walras model solve a system of linear 
equations coming via algebraic manipulations from agent utility maximization 
and supply-demand balance.  In the notation
$\tilde\beta_{ij}=\beta_{ij}/\sumn_{k=0}^n\beta_{ik}$,
this system combines the traditional price normalization equation 
$\sumn_{j=1}^n \tilde p_j =1$ with
	\begin{eqnarray*}
  \tilde x_{ij}-\tilde\beta_{ij}\left(\sum_{j=1}^n\tilde p_j x_{ij}^0\right)&=&
         \tilde\beta_{ij}x_{i0}^0,\; i=1,\ldots,m,\,j=0,\ldots,n\\
     \tilde p_j \sum_{i=1}^m x^0_{ij} -\sum_{k=1}^n\tilde p_k
       \left(\sum_{i=1}^m\tilde\beta_{ij} x_{il}^0\right)
      &=&\sum_{i=1}^m \tilde\beta_{ij}x_{i0}^0,\; j=0,\ldots,n.
	\end{eqnarray*}
The money-denominated prices in Table 2 are derived from this by 
$$
        \bar p^W_j =\frac{\tilde p_j}{\tilde p_0} \text{for} j-1,\ldots,n.
\eqno(5.3)$$
	
{\bf Bigger experiment in the same setting.}  In order to test scalability 
of the trading algorithm, we also applied it to a larger version of 
Example 1 with $n+1=100$ goods and $m=10$ agents.  Despite the much greater 
scope of activity, an equilibrium was again achieved on all runs, which took 
a median of 59345 iterations within 711 seconds of execution time.  
\medskip

The next computational example is aimed differently.  Instead of scaling
upward, we downsize while confronting the algorithm with a challenge.
from another angle.

\proclaim Example 2.
This concerns an economy like the one in Example 1 but with only $n+1=3$ 
goods and $m=3$ agents, and, most importantly, with the initial holdings 
of the agents hugely out of balance.  The Cobb-Douglas exponents and 
initial holdings are indicated in Table 3.  Note that each of the three 
goods is held at first mostly by just one of the three agents.  
	
	\begin{table}[h]  
		\centering
		\begin{tabular}{cc|ccc}
			Agent & Parameter & $j=0$ & $j=1$ & $j=2$ \\
			\hline
			\multirow{2}{*}{$i=1$} & $\beta_1$ & 0.60 & 0.15 & 0.15 \\
			& $x_1^0$ & 10 & 10 & 10 \\
			\hline
			\multirow{2}{*}{$i=2$} & $\beta_2$ & 0.01 & 0.85 & 0.04 \\
			& $x_2^0$ & 2 & 8 & 80 \\
			\hline
			\multirow{2}{*}{$i=3$} & $\beta_3$ & 0.01 & 0.09 & 0.80 \\
			& $x_3^0$ & 2 & 80 & 8
		\end{tabular}
		\caption{Example 2 utility parameters $\beta_{ij}$, 
                  and initial holdings $x^0_{ij}$}
	\end{table}	
	
For this example we report in Table 4 the equilibrium price outcomes of 
two different runs $\nu=1,2$.  The associated final holdings of the agents 
are shown in Table 5.  Despite the initial differences, the results are quite
consistent with each other.  The corresponding prices and holdings coming
from the Walras model for the nonmonetary goods are very different,
however.  This appears all the more to cast doubt on that model's ability 
to reflect genuine market forces.

	\begin{table}[h]  
		\centering
\begin{tabular}{l|rrrrrr}
	 &   $\overline p^{1}$ &  $\overline p^{2}$ &  $\overline p^{W}$ \\
	\hline
	$j=1$ &             0.0475 &             0.0460 &             0.4921 \\
	$j=2$ &             0.0494 &             0.0484 &             0.4614 \\
\end{tabular}
		\caption{Equilibrium prices $\bar{p}^\nu$ for Example 2}
	\end{table}

	\begin{table}[h] 
		\centering
\begin{tabular}{l|rrrrrrrrr}
	 &  $\overline x^{\nu}_{10}$ &  $\overline x^{\nu}_{20}$ 
         &  $\overline x^{\nu}_{30}$ &  $\overline x^{\nu}_{11}$ 
         &  $\overline x^{\nu}_{21}$ &  $\overline x^{\nu}_{31}$ 
         &  $\overline x^{\nu}_{12}$ &  $\overline x^{\nu}_{22}$ 
         &  $\overline x^{\nu}_{32}$ \\
	\hline

	$\nu=1$ &                     13.97 &                      0.01 &                      0.02 &                     73.58 &                     21.33 &                      3.09 &                     70.66 &                      0.96 &                     26.37 \\
	$\nu=2$ &                     13.97 &                      0.01 &                      0.02 &                     75.95 &                     19.08 &                      2.96 &                     72.13 &                      0.85 &                     25.01 \\
	Walrus   &                     13.02 &                      0.48 &                      0.50 &                      6.62 &                     82.23 &                      9.16 &                      7.06 &                      4.13 &                     86.82 \\
\end{tabular}
	\caption{Final holdings $\bar{x}^\nu$ for Example 2}
	\end{table}
	
Our final example has a different character. It employs utility functions
not of Cobb-Douglas type that allow agents to reduce their holdings of 
non-money goods to 0, if they so desire.  As an aid in setting up these 
utility functions we make use of the fact that at any stage of trading the 
holdings $x_{ij}\geq 0$ of the agents satisfy
$$
         \sum_{i=1}^m x_{ij} = s_j, \text{where} 
         s_j=\sum_{i=1}^m x^0_{ij} \text{(total supply of good $j$), hence}
        x_{ij}\in [0,s_j].
\eqno(5.4)$$

\proclaim Example 3.
This concerns an economy with just $m=2$ agents and $n+1=3$ goods.  Each 
agent $i$ has a separable utility function 
$u_i(x_i)=\sum_{j=0}^2 u_{ij}(x_{ij})$ in which the term for the money
good has the form $u_{i0}(x_{i0}) = x_{i0}^{\alpha_i}$ for some 
$\alpha_i\in\left(0,1\right)$.   The terms for the non-money goods $j=1,2$,
have instead the form
$$
	u_{ij}(x_{ij})=a_{ij}x_{ij}-\frac{1}{2}b_{ij}x_{ij}^2
        \text{for}  x_{ij}\in[0,s_j], \text{with}
      a_{ij}>0,\; b_{ij}>0, \; \frac{a_{ij}}{b_{ij}}>s_j,
\eqno(5.5)$$
where $s_j$ is the total supply of good $j$ as in (5.4).\footnote{
   The conditions in (5.5) ensure that $u_{ij}$ is a strongly concave 
   increasing quadratic function on the interval $[0,s_j]$ and has positive 
   slope still at $s_j$ itself.  It is mathematically possible then to 
   build an extension beyond $s_j$ that is strongly concave, increasing and 
   continuously twice differentiable, so as to end up with a utility 
   function $u_i$ that fits our general assumptions.  But because of (5.4),
   that extension can be left implicit.}
The values of the utility parameters are indicated in Table 6.  Note that
for good $j=2$ a shift is introduced between the two agents of an amount 
$\mu^\nu$ that depends on the run index $\nu$,
	
	\begin{table}[H] 
		\centering
\begin{tabular}{l|rrrr}
	&  $ p^{0,\nu}_{i1}$ &  $\overline p^{\nu}_{i1}$ &  $ p^{0,\nu}_{i2}$ &  $\overline p^{\nu}_{i2}$ \\
	\hline
	$\nu=0,i=1$ &            31.4643 &                   18.9509 &            21.3957 &                   18.1378 \\
	$\nu=0,i=2$ &             2.5298 &                   18.9507 &             3.2888 &                   18.1379 \\
	$\nu=1,i=1$ &            31.4643 &                   17.4484 &            28.9471 &                   15.7335 \\
	$\nu=1,i=2$ &             2.5298 &                   21.0987 &             1.7709 &                   15.7335 \\
\end{tabular}
	\caption{ Utility parameters and initial holdings in 
                   Example 3, where $\mu^\nu=3\nu$}
       \end{table}

	\begin{table}[]  
		\centering
		\begin{tabular}{cc|cc|cc}
			Trial & Agent &\multicolumn{2}{c}{Price threshold good 1}&\multicolumn{2}{|c}{Price threshold good 2}\\
			&&{$p_{i1}^{0,\nu}$} & ${\bar p}_{i1}^{\nu}$ & {$p_{i2}^{0,\nu}$} & {${\bar p}_{i2}^{\nu}$} \\
		\hline
		$\nu=0$ & $i=1$ & 31.4643 & 18.9509 & 21.3957 & 18.1378 \\
		$\nu=0$ & $i=2$ & 2.5298 & 18.9507 & 3.2888 & 18.1379 \\
		\hline
		$\nu=1$ & $i=1$ & 31.4643 & 19.6394 & 25.1714 & 17.5170 \\
		$\nu=1$ & $i=2$ & 2.5298 & 19.6393 & 2.5298 & 17.5172 \\
		\hline
		$\nu=2$ & $i=1$ & 31.4643 & 17.4484 & 28.9471 & 15.7335 \\
		$\nu=2$ & $i=2$ & 2.5298 & 21.0987 & 1.7709 & 15.7335
		\end{tabular}
	\caption{Initial and final prices for Example 3}
	\end{table}

Note that this family of utility functions has marginal utility 
going to $\infty$ as money goes to 0, but not as other goods go to 0.  
The sets of admissible holdings $X_i$ in (2.3) come out therefore as
$$
 X_i =\lset x_i \mset x_{i0}>0,\, x_{i1} \geq 0,\, x_{i2}\geq 0 \rset.
\eqno(5.6)$$
This puts us outside the domain of the guaranteed convergence of the
trading algorithm in Theorem 3, which depended on assumption (4.1), but we
wanted to see anyway what might happen computationally in such circumstances.

Following the rules in the previous sections, the price thresholds  
$p_{ij}(x_i)$ in (2.8) can be explicitly computed as
$$
		p_{ij}(x_{i})= 
  \frac{a_{ij}-b_{ij}x_{ij}}{\alpha_{i}x_{i0}^{\alpha_i-1}},\quad j=1,2
\eqno(5.7)$$
we can also compute the solution to the seller's problem according to (2.14),
$$\eqalign{
   \xi_j^{+}(x_i,\pi_j)=\argmax_{\xi_j}\!\lset u_i(x_i+\xi_j[\pi_j,-1])
       \mset x_i+\xi_j[\pi_j,-1]\in X_i\rset  &\cr
 \hskip52pt =\displaystyle{ \argmax_{0\leq\xi_j\leq x_{ij}}}
       \lset (x_{i0}+\xi_j\pi_j)^{\alpha_i}
     +a_{ij}(x_{ij}-\xi_j)-\frac{1}{2}b_{ij}(x_{ij}-\xi_j)^2 \rset &}
\eqno(5.8)$$
and, correspondingly the buyer's problem from (2.20)
$$\eqalign{
 \xi_j^{-}(x_i,\pi_j)=\argmax_{\xi_j}\!\lset u_i(x_i-\xi_j[\pi_j,-1])
     \mset\,x_i-\xi_j[\pi_j,-1]\in X_i\rset &\cr
 \hskip52pt  =\displaystyle{\argmax_{0\leq\xi_j < \frac{x_{i0}}{\pi_j}}}
   \lset (x_{i0}-\xi_j\pi_j)^{\alpha_i}
  +a_{ij}(x_{ij}+\xi_j)-\frac{1}{2}b_{ij}(x_{ij}+\xi_j)^2 \rset. &}
\eqno(5.9)$$

We ran three instances of this example in which agent 1 starts with none of
good $j=1$, agent 2 starts with almost no money, and the initial amount of
good $j=2$ shifts from agent 1 to agent 2 in dependence on the run index,
$\nu=0,\,1,\,2$.  The holdings and prices on termination are reported 
in Table 7.

Two different behaviors are seen.  On the one hand, in runs $\nu=0,\,1$, 
our algorithm converged to an equilibrium within an average of 174 iterations 
and 0.512 seconds, despite initial holdings not all being positive and 
utility functions not satisfying assumption (4.1).  However, for trial 
$\nu= 2$, we observed a lack of convergence of the price thresholds for 
good $j=1$, as can be seen in Figure 1; this is the good which agent 1 had 
none of at the start.  {\it The failure to reach an equilibrium in this 
instance indicates that, for economies with admissible holdings not 
necessarily positive, as allowed by their utility functions not fitting 
with assumption (4.1), something more needs to be added to the trading 
scheme to avoid stagnation.}
	
\begin{figure}[H]
	\includegraphics[width=1.0\textwidth]{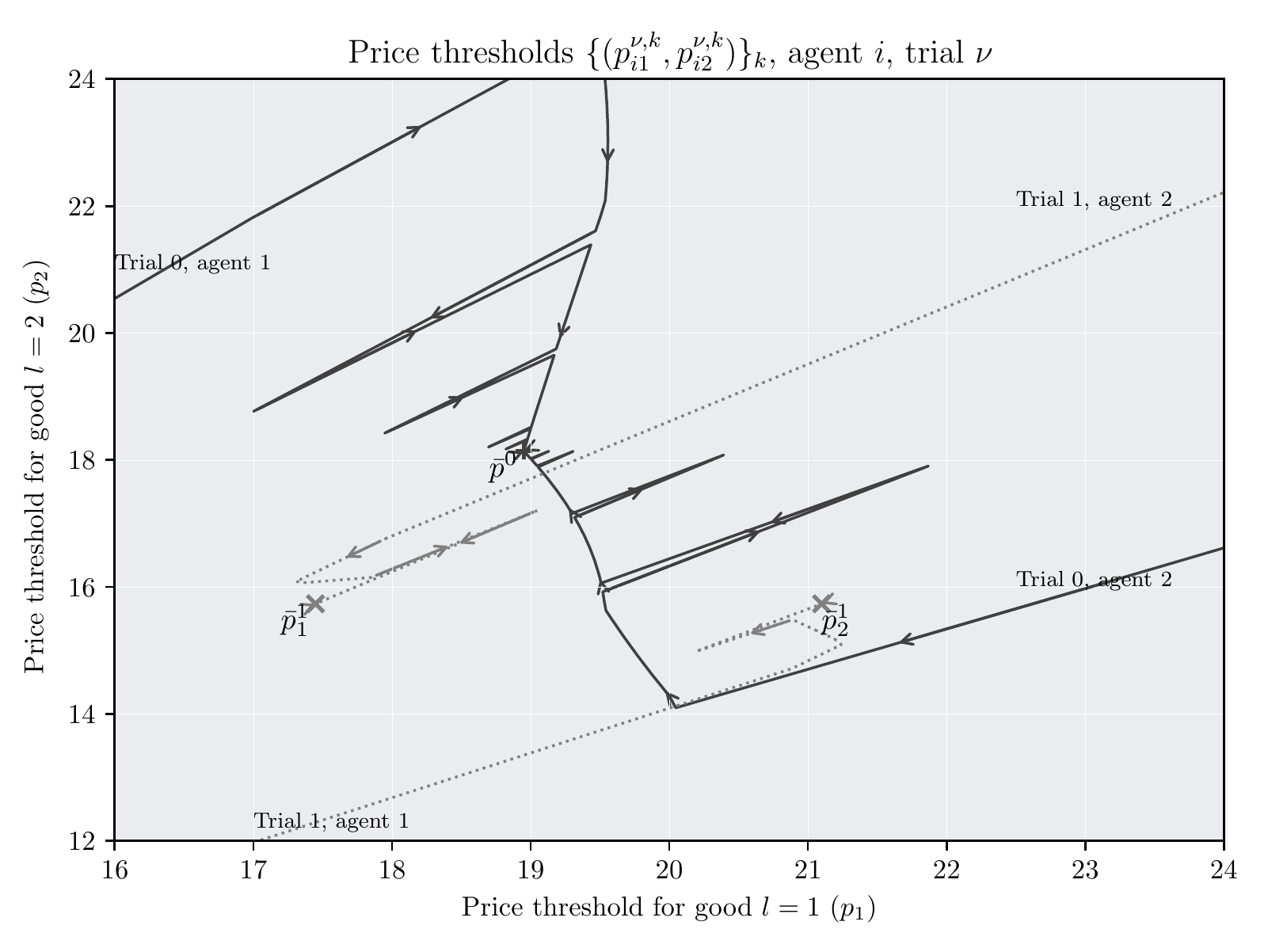}
  \caption{Trial 0.  Successful example, despite violation of assumption
   (4.1).  Equilibrium reached as the price threshold curve for agent 1,  
   $\{p_1^{0,k}\}$, met the curve for agent 2, $\{p_2^{0,k}\}$.  Trial 1.  
   Failed example in the absence of (4.1).  Both curves remained far from 
   each other even though price thresholds for good 2 came to coincide, 
       ($\bar{p}_{1,2}^{0}=\bar{p}_{2,2}^{1}$).}
\end{figure}
	
The bilateral trade dynamics between the agents is depicted below in 
Figure 2, where we plot the Edgeworth box for the nonmonetary pair of 
goods in the economy.    
	
\begin{figure}[ht!]
	\includegraphics[width=1.0\textwidth]{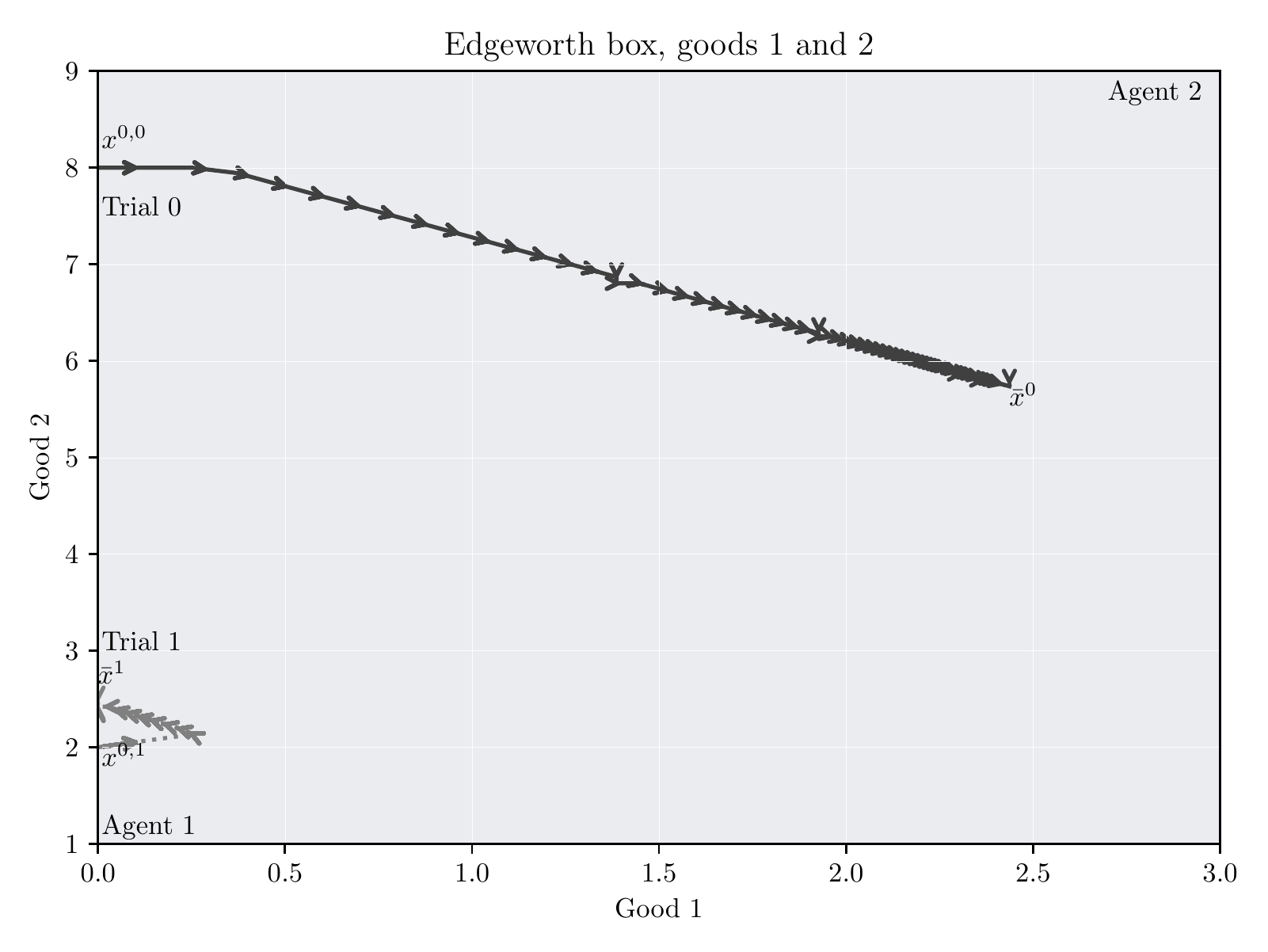}
 \caption{Trial 0.  Starting from an initial endowment in the upper left 
  region of the box, the agents bilaterally trade until they reach 
  equilibrium allocations ($\bar{x}^0$).  Trial 1.  Agents enter the economy 
  with allocations in the lower left corner of the box, but bilateral 
  exchanges lead them to a boundary point, which does not correspond to an
  equilibrium.}
\end{figure}

 

\newpage


\begin{thebibliography}{99} 

%
\bibitem{Balasko0} {\sc Y.\ Balasko,}
Some results on uniqueness and on stability of equilibrium in general
equilibrium theory.
{\sl J.\ Mathematical Economics\/} 2 (1975), 95--118.

%
\bibitem{Bottazzi} {\sc J.-M.\ Bottazzi,}
Accessibility of Pareto optima by Walrasian exchange processes.
{\sl Journal of Mathematical Economics\/} 23 (1994), 583--603.

%
\bibitem{Debreu-Value} {\sc G.\ Debreu,}
{\sl Theory of Value:  An Axiomatic Analysis of Economic Equilibrium.}
Wiley, New York, 1959.

%
\bibitem{DubeyShubik} {\sc P.\ Dubey, M.\ Shubik,}
A strategic market game with price and quantity strategies.
{\sl Journal of Economics (Zeitschrift f\"ur National\"okonomie)\/} 40
(1980), 25--34.

%
\bibitem{DubeySahiShubik} {\sc P.\ Dubey, S.\ Sahi, M.\ Shubik,}
Money as minimal complexity.
{\sl Games and Economic Behavior\/} 108 (2018), 432--451.
%
\bibitem{Eckalbar} {\sc J.\ C.\ Eckalbar,}
Bilateral trade in a monetized pure exchange economy.  {\sl Economic
Modeling\/} 3 (1986), 135--139.

%
\bibitem{Feldman} {\sc A.\ M.\ Feldman, }
Bilateral trading processes, pairwise optimality, and Pareto optimality.
{\sl The Review of Economic Studies\/} 46 (1973), 463--473.

%
\bibitem{Flaam1} {\sc S.\  Fl{\aa}m, }
Bilateral exchange and competitive equilibrium.
{\sl Set-valued and Variational Analysis\/} 24 (2016), 1--11.

%
\bibitem{Flaam2} {\sc S.\  Fl{\aa}m, K.\ Gramstad,}
Reaching Cournot-Walras equilibrium. {\sl ESAIM: Proceedings and Surveys\/}
57 (2017), 12--22.

%
\bibitem{Flaam3} {\sc S.\  Fl{\aa}m, }
Emergence of price-taking behavior.  
{\sl Economic Theory} {\bf 70} (2020), 847--870.

%
\bibitem{Flaam4} {\sc S.\  Fl{\aa}m, }
Towards competitive equilibrium by double auctions. {\sl Pure and Applied
Functional Analysis\/} {\bf 6} (2021), 1211--1225.

%
\bibitem{Flaam5} {\sc S.\  Fl{\aa}m, }
Market equilibria and money.  {\sl Fixed Point Theory Algorithms Sci.\
Eng.}, doi.org/10.1186/s13663-021-00705-4.

%
\bibitem{Gintis} {\sc S.\  H.\ Gintis, }
The dynamics of general equilibrium.  {\sl The Economic Journal\/} 
{\bf 117} (2007), 1280--1307.

%
\bibitem{GEI} {\sc A.\ Jofr\'e, R.\ T.\ Rockafellar, R.\ J-B Wets, }
General equilibrium with financial markets and retainability.
{\sl Economic Theory\/} 63 (2017), 309--345.

%
\bibitem{Stability} {\sc A.\ Jofr\'e, R.\ T.\ Rockafellar, R.\ J-B Wets, }
On the stability and evolution of economic equilibrium.  
{\sl Convex Analysis\/} 30 (2023).  Downloadable as \#246 from
sites.math.washington.edu/$\sim$rtr/mypage.html.   

%
\bibitem{Keisler} {\sc H.\ J.\ Keisler,} 
Getting to a competitive equilibrium.  
{\sl Econometrica\/} 64 (1996), 29--49.

%
\bibitem{Kitti} {\sc M.\ Kitti,}
Convergence of iterative t\^atonnement without price normalization.
{\sl Journal of Economic Dynamics and Control\/} 34 (2011), 1077--1081.

%
\bibitem{Levando} {\sc D.\ Levando,}
A Survey of Strategic Market Games.
{\sl Economic Annals\/} LVII, no.\ 194 (2012).

%
\bibitem{Mandel} {\sc A.\ Mandel, H.\ Gintis,}
Decentralized pricing and strategic stability of Walrasian geneeral
equilibrium.  {\sl J.\ Math.\ Econ.} {\bf 63} (2016), 84--92.

%
\bibitem{Ostroy} {\sc J.\ M.\ Ostroy,}
The informational efficiency of monetary exchange.
{\sl The American Economic Review\/} 63 (1973), 597--610.

%
\bibitem{WalrasExchanges} {\sc R.\ T.\ Rockafellar,}
Optimization and decentralization in the mathematics of economic
equilibrium.  {\sl Proc.\ of the Internat'l Conf.\ on
Nonlinear Analysis and Convex Analysis \& Internat'l Conf.\ on
Optimization Techniques and Applications\/} 11 (Yokohama Publishers, 2021), 
199-211.  Downloadable as \#255 from
sites.math.washington.edu/$\sim$rtr/mypage.html.   

%
\bibitem{Preferences} {\sc R.\ T.\ Rockafellar,}
Variational analysis of preference relations and their utility
representations.  {\sl Pure and Applied Functional Analysis}, accepted.
Downloadable as \#263 from
sites.math.washington.edu/$\sim$rtr/mypage.html.   

%
\bibitem{Rasooly} {\sc I.\ Rasooly,}
Competitive equlibrium and the double auction.   Preprint 2022.  

%
\bibitem{Saari} {\sc D.\ G.\ Saari,}
Iterative price mechanisms.
{\sl Econometrica\/} 53 (1985), 1117--1132,

%
\bibitem{Sorin} {\sc S.\ Sorin,} 
Strategic games with exchange rates.
{\sl J.\ Economic Theory\/} 69 (1996), 431--446.

%
\bibitem{ShapleyShubik} {\sc L.\ Shapley, M.\ Shubik,}
Trade using one commodity as a means of payment.
{\sl J.\ Political Economy\/} 85 (1977), 937--968.

%
\bibitem{Starr} {\sc R.\ M.\ Starr,}
The structure of exchange in barter and monetary economies.
{\sl The Quarterly Journal of Economics\/} 86 (1972), 290--302.

%
\bibitem{Uzawa} {\sc H.\ Uzawa,}
Walras t\^atonnement in the theory of exchange.
Technical report no.\ 69, Project on Efficiency in Decision Making in
Economic Systems, Stanford University, 1959.

%
\bibitem{Walker} {\sc D.\ A.\ Walker,}
Walras's theories of t\^atonnement.
{\sl Journal of Political Economy\/} 95 (1987), 758--774. 

\end{thebibliography}
\end{document}